\documentclass[a4paper,12pt,final]{amsart}
\usepackage{times,a4wide,mathrsfs,amssymb,amsmath,amsthm,dsfont}

\newcommand{\C}{\mathbb{C}}
\newcommand{\ZZ}{\mathbb{Z}}

\newcommand{\QQ}{\mathbb{Q}}
\newcommand{\NN}{\mathbb{N}}
\newcommand{\PP}{\mathbb{P}}

\newcommand{\OO}{\mathcal O}

\newcommand{\Sy}{\mathfrak S}

\newcommand{\MM}{\mathcal M}

\newcommand{\rom}{\romannumeral}

    \newcommand{\one}{\mathds{1}}

\DeclareMathOperator{\sym}{Sym}
\DeclareMathOperator{\jac}{Jac}

\newtheorem{theorem}{Theorem}[section]

\newtheorem{lemma}[theorem]{Lemma}

\newtheorem{corollary}[theorem]{Corollary}

\newtheorem{conjecture}[theorem]{Conjecture}

\theoremstyle{definition}

\newtheorem{convention}{Conventions}

\newtheorem{nonumberingt}{Acknowledgements}

\newtheorem{nonumberingcon}{Conflict of interest statement}

\begin{document}
\author[Robert Laterveer]
{Robert Laterveer}

\address{Institut de Recherche Math\'ematique Avanc\'ee,
CNRS -- Universit\'e 
de Strasbourg,\
7 Rue Ren\'e Des\-car\-tes, 67084 Strasbourg CEDEX,
FRANCE.}
\email{robert.laterveer@math.unistra.fr}

\title[Zero-cycles on self-products of varieties]{Zero-cycles on self-products of varieties: some elementary examples verifying Voisin's conjecture}

\begin{abstract} An old conjecture of Voisin describes how zero-cycles on a variety $X$ should behave when pulled-back to the self-product 
$X^m$ for $m$ larger than the geometric genus of $X$.
Using complete intersections of quadrics, we give examples of varieties in any dimension and with arbitrarily high geometric genus that verify Voisin's conjecture.
\end{abstract}

\keywords{Algebraic cycles, Chow groups, motives, Voisin conjecture, complete intersections of quadrics}
\subjclass[2010]{Primary 14C15, 14C25, 14C30.}

\maketitle

\section{Introduction}

Given a smooth projective variety $X$ over $\C$, let $A^i(X)_{\ZZ}:=CH^i(X)_{\ZZ}$ denote the Chow groups of $X$ (i.e. the groups of codimension $i$ algebraic cycles on $X$ with $\ZZ$-coefficients, modulo rational equivalence \cite{F}). Let $A^i_{hom}(X)_{\ZZ}$ denote the subgroup of homologically trivial cycles.     

The Bloch--Beilinson--Murre conjectures form a kind of Rosetta Stone, allowing to translate cohomological statements into conjectural Chow-theoretic statements \cite{J2}, \cite{J4}, \cite{Mur}, \cite{Kim}, \cite{MNP}, \cite{Vo}. The following particular instance of such a translation was first formulated by Voisin:

\begin{conjecture}[Voisin 1993 \cite{V9}]\label{conj} Let $X$ be a smooth projective $n$-dimensional variety with $H^j(X,\OO_X)=0$ for $0<j< n$.
Let $m$ be an integer strictly larger than the geometric genus $p_g(X):=\dim H^n(X,\OO_X)$. Then for any zero-cycles $a_1,\ldots,a_m\in A^n_{hom}(X)_{\ZZ}$, there is an equality
  \[ \sum_{\sigma\in\Sy_m} \bigl(\hbox{sgn}(\sigma)\bigr)^{n+1}\,  a_{\sigma(1)}\times\cdots\times a_{\sigma(m)}=0\ \ \ \hbox{in}\ A^{mn}(X^m)_{\ZZ}\ .\]
  (Here $\Sy_m$ is the symmetric group on $m$ elements, and $ \hbox{sgn}(\sigma)$ is the sign of the permutation $\sigma$.
  The notation $a_1\times\cdots\times a_m$ is shorthand for the zero-cycle $p_1^\ast(a_1)\cdot p_2^\ast(a_2)\cdots p_m^\ast(a_m)$ on 
  $X^m$, where the $p_j\colon X^m\to X$ are the various projections.)
  \end{conjecture}
  
  This conjecture is the translation in Chow-language of the fact that the Hodge structure 
  \[ \wedge^m H^n(X,\QQ)\ \ \subset\ H^{mn}(X^m,\QQ)\] 
  has Hodge coniveau $>0$ for $m>p_g(X)$ (see \cite{V9} or \cite[Section 4.3.5.2]{Vo} for more detailed motivation).
  
In case $p_g(X)=0$, Conjecture \ref{conj} predicts that $A^n(X)_\ZZ\cong\ZZ$ (this is a form of Bloch's conjecture \cite{B}). In case $p_g(X)=1$ (e.g., $X$ is a Calabi--Yau variety), the conjecture takes on a particularly appealing form: it predicts that any 2 degree 0 zero-cycles $a,a^\prime\in A^n_{hom}(X)_{\ZZ}$ verify
  \[  a\times a^\prime = (-1)^n \, a^\prime\times a\ \ \ \hbox{in}\ A^{2n}(X\times X)_{\ZZ}\ .\]
This
is still open for a general K3 surface (in fact, I am not aware of a K3 surface with Picard number $<9$ which is known to verify Voisin's conjecture).
Examples of surfaces of geometric genus $1$ verifying the conjecture (or a variant conjecture) are given in \cite{V9}, \cite{16.5}, \cite{19}, \cite{21}, \cite{Zan}. Examples
of other varieties verifying the conjecture are given in \cite{V9}, \cite{32}, \cite{17}, \cite{24.4}, \cite{24.5}, \cite{CF}, \cite{BLP}, \cite{LV}, \cite{Ch}, \cite{Bur}.

 The main result of this note is that Voisin's conjecture is true for certain complete intersections of quadrics:

\begin{theorem}\label{main} Given $g,r\in\NN$ and distinct numbers $\lambda_0,\ldots,\lambda_{2g+1}\in\C$,
let $X\subset\PP^{2g+1}(\C)$ be the complete intersection defined by 
    \[   \begin{cases}  x_0^2+x_1^2+\cdots \  \cdots \ \cdots+ x_{2g+1}^2&=0\ ,\\
                          \lambda_0 x_0^2+\lambda_1 x_1^2+\cdots + \lambda_{2g+1}x_{2g+1}^2&=0\ ,\\
                            \lambda_0^2 x_0^2+\lambda_1^2 x_1^2+\cdots + \lambda_{2g+1}^2 x_{2g+1}^2&=0\ ,\\
                          \ \ \ \vdots\\
                         \ \ \  \vdots\\
                           \lambda_0^r x_0^2+\lambda_1^r x_1^2+\cdots + \lambda_{2g+1}^r x_{2g+1}^2&=0\ .\\     
                           \end{cases}\]
                           Then $X$ is a smooth projective variety, and Conjecture \ref{conj} is true for $X$.
\end{theorem}

These complete intersections are admittedly very special (and unfortunately the argument proving Theorem \ref{main} does {\em not\/} apply to general complete intersections of quadrics), but at least they provide examples of any dimension, and with $p_g$ arbitrarily high, verifying Voisin's conjecture. The complete intersections of Theorem \ref{main} are rather similar to the Calabi--Yau varieties of \cite{LV}: both are related to products of curves, and hence to abelian varieties. Thus, to prove Theorem \ref{main} we can reduce to a problem concerning zero-cycles on abelian varieties; this last problem can be solved thanks to recent work of Vial \cite{Ch}.

 As a consequence of Theorem \ref{main}, certain instances of the generalized Hodge conjecture are verified:
 
 \begin{corollary}\label{ghc} Let $X$ be as in Theorem \ref{main}, and $m>p_g(X)$. The sub-Hodge structure
   \[ \wedge^m H^{n}(X,\QQ)\ \ \subset\ H^{mn}(X^m,\QQ) \]
   is supported on a divisor.
   \end{corollary}

   \vskip0.6cm

\begin{convention} In this note, the word {\sl variety\/} will refer to a reduced irreducible scheme of finite type over $\C$. A {\sl subvariety\/} is a (possibly reducible) reduced subscheme which is equidimensional. 

{\bf Unless indicated otherwise, all Chow groups will be with rational coefficients}: we will denote by $A_j(X)$ (resp. $A_j(X)_{\ZZ}$) the Chow group of $j$-dimensional cycles on $X$ with $\QQ$-coefficients (resp. $\ZZ$-coefficients); for $X$ smooth of dimension $n$ the notations $A_j(X)$ and $A^{n-j}(X)$ are used interchangeably. 
The notation $A^j_{hom}(X)$ will be used to indicate the subgroups of homologically trivial cycles.

We will write $\MM_{\rm rat}$ (and $\MM_{\rm hom}$) for the contravariant category of pure motives with respect to rational equivalence (resp. homological equivalence), as in \cite{Sc}, \cite{MNP}.
 \end{convention}

 \section{The proof}
 
     
   This section contains the proof of the main result:

    \begin{proof}(of Theorem \ref{main}.) This is based on results of Terasoma, who has made an in-depth study of this kind of complete intersection \cite[Sections 2.4, 2.5 and 2.6]{Tera} (NB: in loc. cit., the variety $X$ is denoted $X_{r,g}$ which is (by definition) the same as $X_{r,2g+2,2}$). 
    
    First, the non-singularity of $X$ is readily checked using the Vandermonde determinant (cf. \cite[Proposition 2.4.1]{Tera}).
    
    Next, to verify Conjecture \ref{conj} for $X$, we need to understand the Chow group of zero-cycles $A^n(X)_\ZZ$ (where we write $n:=2g-r$ for the dimension of $X$). Thanks to a theorem of Roitman \cite{Ro}, it suffices to prove that Conjecture \ref{conj} is true for $A^n(X)$, the Chow group with $\QQ$-coefficients.
  From Terasoma's work \cite[Corollary 2.5.3]{Tera}, we know there exist a finite number of hyperelliptic curves $C_\chi$ and a correspondence $\Gamma$ inducing an isomorphism in cohomology
     \begin{equation}\label{ter} \Gamma_\ast\colon\ \ H^n_{prim}(X,\QQ)\ \xrightarrow{\cong}\ \bigoplus_\chi \wedge^n H^1(C_\chi,\QQ)\ .\end{equation}
  Because curves and complete intersections verify the Lefschetz standard conjecture, we know that the inverse is also induced by a correspondence. (This is well-known, cf. for instance \cite[Proof of Proposition 1.1]{V4}, where I first learned this. In a nutshell, the argument is as follows: the Lefschetz standard conjecture for $X$ implies that there is a correspondence $s\in A^n(X\times X)$ inducing a polarization on $H^n_{prim}(X,\QQ)$. A Hodge-theoretical argument \cite[Lemma 1.6]{V4} then gives that $\alpha:=\Gamma\circ s \circ {}^t \Gamma$ induces an automorphism of $\bigoplus_\chi \wedge^n H^1(C_\chi,\QQ)$. By the Cayley--Hamilton theorem, the inverse to $\alpha_\ast$ is given by $P(\alpha)_\ast$ for some rational polynomial $P$. Then $\Gamma\circ s\circ {}^t \Gamma\circ P(\alpha)$ acts as the identity on $\bigoplus_\chi \wedge^n H^1(C_\chi,\QQ)$, and so the correspondence $s\circ {}^t \Gamma\circ P(\alpha)$ gives an inverse to $\Gamma_\ast$. )
     
     This means that there is an isomorphism of homological motives
     \begin{equation}\label{homo} [\Gamma]\colon\ \ h^n_{prim}(X)\ \xrightarrow{\cong}\ \bigoplus_\chi \sym^n h^1(C_\chi)\ \ \ \hbox{in}\ \MM_{\rm hom}\ .\end{equation}
     Let $h^n_{prim}(X)\in\MM_{\rm rat}$ be the Chow motive such that 
       \[ h(X)=h^n_{prim}(X)\oplus \bigoplus\one(\ast)\ \ \ \hbox{in}\ \MM_{\rm rat}\] 
   (that is, $h^n_{prim}(X)$ is defined by the projector
     $\pi^{n,prim}_X:=\Delta_X -{1\over 2^{r+1}} \sum_j h^j\times h^{n-j}$, where $h\in A^1(X)$ denotes a hyperplane section). The variety $X$ is isomorphic to a quotient $D^{2g-r}/G$, where $D$ is a smooth curve \cite[Theorem 2.4.2]{Tera}, and so $X$ is Kimura finite-dimensional \cite{Kim}. We recall Kimura's fundamental nilpotence theorem \cite{Kim}, which implies that the functor from
   finite-dimensional Chow motives to homological motives is {\em conservative\/}, i.e. detects isomorphisms.
     The isomorphism \eqref{homo} can thus be upgraded to an isomorphism of Chow motives
     \begin{equation}\label{chow}  h^n_{prim}(X)\ \xrightarrow{\cong}\ \bigoplus_\chi \sym^n h^1(C_\chi)\ \ \ \hbox{in}\ \MM_{\rm rat}\ .     \end{equation}
 Let $J_\chi:=\jac(C_\chi)$ denote the Jacobian of the curve $C_\chi$. There are isomorphisms $h^1(C_\chi)\cong h^1(J_\chi)$. Moreover, $\sym^n h^1(B)\cong h^n(B)$ for any abelian variety $B$ (here $h^\ast(B)$ refers to the Deninger--Murre decomposition of abelian schemes \cite{DM}), and so the isomorphism \eqref{chow} induces an isomorphism
   \[  h^n_{prim}(X)\ \xrightarrow{\cong}\ \bigoplus_\chi h^n(J_\chi)\ \ \ \hbox{in}\ \MM_{\rm rat}\ .\]     
   
 The following lemma (which we state separately for possible future reference) now closes the proof of the theorem:
 
 \begin{lemma} Let $X$ be a smooth projective variety of dimension $n$, and assume that
   \begin{equation}\label{chowlem} h(X)\cong \bigoplus_\chi h^n(B_\chi)\oplus \bigoplus \one(\ast)\ \ \ \hbox{in}\ \MM_{\rm rat}\ ,\end{equation}
   where $B_\chi$ are abelian varieties, and $h^\ast(B_\chi)$ refers to the Deninger--Murre decomposition of the motive of abelian schemes \cite{DM}. Then Conjecture \ref{conj} is true for $X$.
   \end{lemma}  
   
       
    

It remains to prove the lemma. (NB: In proving the lemma, we use some results concerning the Chow motive of an abelian variety $B$ and the Beauville decomposition $A^\ast_{(\ast)}(B)$ of the Chow ring. For a quick survey containing all the results we use, one could look at \cite[Section 5]{Sc}.)

    Taking Chow groups on both sides of \eqref{chowlem}, we obtain an isomorphism
      \begin{equation}\label{here} A^n_{hom}(X)\ \xrightarrow{\cong}\ \bigoplus_\chi  A^n\bigl(    h^n(B_\chi)\bigr) \ .\end{equation}
      Writing $g_\chi:=\dim B_\chi$, let us
      proceed to analyze the various summands in \eqref{here}:
      
      \begin{itemize}
      
      \item In case $n>g_\chi$, the summand $A^n\bigl(    h^n(B_\chi)\bigr)$ is zero for dimension reasons. The piece $A^{g_\chi}_{(n)}(B_\chi)$ of the Beauville decomposition is also zero.
      
      \item In case $n=g_\chi$, we have that $A^n\bigl(    h^n(B_\chi)\bigr)= A^{g_\chi}_{(n)}(B_\chi)$, since the Deninger--Murre decomposition $h^\ast(B)$ induces the Beauville decomposition $A^\ast_{(\ast)}(B)$ in the sense that $A^j_{(s)}(B)=A^j(h^{2j-s}(B))$ for all integers $j$ and $s$ and any abelian variety $B$.
      
      \item Finally, let us assume $n<g_\chi$, say $n=g_\chi-s$ where $s>0$. In this case, the work of K\"unnemann \cite{Kun} provides a ``hard Lefschetz''  isomorphism of motives
      \[ h^n(B_\chi)\ \xrightarrow{\cong}\ h^{n+2s}(B_\chi)(s)\ \ \ \hbox{in}\ \MM_{\rm rat}\ .\]
      Taking Chow groups, this induces an isomorphism
       \[  A^n\bigl(    h^n(B_\chi)\bigr)      \ \xrightarrow{\cong}\  A^{n+s}\bigl( h^{n+2s}(B_\chi)\bigr)=A^{n+s}_{(n)}(B_\chi)=A^{g_\chi}_{(n)}(B_\chi)\ .\]
       \end{itemize}
       
   In view of this analysis, we can rewrite \eqref{here} in order to have zero-cycles on both sides:
      \begin{equation}\label{here2} A^n_{hom}(X)   \ \xrightarrow{\cong}\ \bigoplus_\chi     A^{g_\chi}_{(n)}(B_\chi)\ .\end{equation}
      
  We now invoke a result of Vial:
   
   \begin{theorem}[Vial \cite{Ch}]\label{via} Let $B$ be an abelian variety of dimension $g$, and let $\pi^j_B$ denote the Deninger--Murre projectors \cite{DM}.
   
   \noindent
   (\rom1) Assume $n\in\NN$ is even. Then
       \[  (\wedge^m \, \pi_B^{2g-n})_\ast A_0(B^m)=0\ \ \ \forall\ m>{g\choose n}\ .\]
   
   \noindent
   (\rom2) Assume $n\in\NN$ is odd. Then
        \[  (\sym^m \, \pi_B^{2g-n})_\ast A_0(B^m)=0\ \ \ \forall\ m>{g\choose n}\ .\]   
   \end{theorem} 
   
   (This is \cite[Theorem 4.1]{Ch}. Vial's result generalizes a result of Voisin \cite[Example 4.40]{Vo}, which was the case $n=g$ of Theorem \ref{via}.)

                    Armed with Vial's result, we are in position to prove the lemma. Let us treat in detail the case $n$ even (the case $n$ odd is similar and will be left to the zealous reader). 
                    Let us write $\pi^{prim}_X$ for the projector such that $(X,\pi^{prim}_X,0)\cong \oplus_\chi h^n(B_\chi)$.
                     Given an integer $m>p_g(X)$, using \eqref{here2} we find 
    \begin{equation}\label{this}  (\wedge^m \, \pi^{prim}_X)_\ast   A_0(X^m) \   \cong\  \bigoplus_{\sum m_\chi=m}  \bigotimes_\chi (\wedge^{m_\chi} \, \pi_{B_\chi}^{2g_\chi-n})_\ast    A_0\bigl( B_\chi^{m_\chi})   \  . \end{equation}
  We note that \eqref{chowlem} gives us an equality
    \[  p_g(X) ={\displaystyle\sum_\chi} {g_\chi\choose n}\ .\]
    Since (by assumption) $m>p_g(X)$, it follows that in each partition $m=\sum m_\chi$
    there exists some $\chi$, say $\chi_0$, such that 
      \[ m_{\chi_0}>    {g_{\chi_0}\choose n}\ .\]

   Theorem \ref{via}(\rom1) applied to $B_{\chi_0}$ guarantees that in the sum \eqref{this} each summand  
       \[  \bigotimes_\chi (\wedge^{m_\chi}\, \pi_{B_\chi}^{2g_\chi-n})_\ast    A_0\bigl( B_\chi^{m_\chi})\] 
     vanishes. It follows that the whole sum \eqref{this} vanishes, i.e.
   \[      (\wedge^m \, \pi^{prim}_X)_\ast   A_0(X^m)=0 \ .\]
   Since $(\pi^{prim}_X)_\ast A_0(X)=A^n_{hom}(X)$, this proves Conjecture \ref{conj} for $X$, in case $n$ is even.
   
   The argument for $n$ odd is similar, the difference being that one considers $\sym^m \pi^{prim}_X$ instead of $\wedge^m \pi^{prim}_X$, and one relies on part (\rom2) of Theorem \ref{via} instead of part (\rom1).
              \end{proof}

   
   \begin{proof}(of Corollary \ref{ghc}.) As Voisin had already remarked \cite[Corollary 3.5.1]{V9}, this is implied by the truth of Conjecture \ref{conj} for $X$ (the implication can be seen using the Bloch--Srinivas argument \cite{BS}; this is explained in detail in \cite[Corollary 2.7]{32}).
   \end{proof}

\vskip0.5cm
\begin{nonumberingt} Thanks to Kai for enjoyable bike trips in the Alsace countryside. Thanks to the referee for many constructive comments that helped to improve the presentation.
\end{nonumberingt}

\vskip0.5cm
\begin{nonumberingcon} The author states that there is no conflict of interest.
\end{nonumberingcon}

\end{document}